\documentclass[11pt]{amsart}
\newcommand{\nsa}{deterministic nested stack automaton with limited erasing}
\newcommand{\A}{A}                      % automaton 
\newcommand{\CG}{\mathcal A}          % configuration graph
\newcommand{\G}{G}                      % group
\newcommand{\CD}{\mathcal G}          % Cayley diagram of a group  
\newcommand{\T}{T}                      % memory tree
\newcommand{\mtrees}{\mathfrak T}     % class of all memory trees
\newcommand{\Mnsa}{M_{\text{nsa}}}    % monoid of memory operations
\newcommand{\e}{\epsilon}               % The empty word 
\newcommand{\de}{\delta} 
\newcommand{\De}{\Delta}
\newcommand{\DD}{\Delta^*} 

\newcommand{\si}{\sigma} 
\newcommand{\Si}{\Sigma} 
\newcommand{\SSi}{\Sigma ^*} 
\newcommand{\inv}{^{-1}} 
\newcommand{\ovr}[1]{\overline #1} 
\newcommand{\cents}{\vert\!\! c }       % ugly cents sign
\newcommand{\df}[1]{{\it #1}}

\newtheorem{theorem}{Theorem}[section] 
\newtheorem{lemma}[theorem]{Lemma}

\theoremstyle{definition} 
\newtheorem{definition}[theorem]{Definition}
\theoremstyle{remark}

\let\phi\varphi

\title[NESTED STACK AUTOMATON GROUPS] 
{On groups whose word problem is solved by a nested stack automaton}

\author[Gilman]{Robert Gilman}

\thanks{The first author was partially supported by National Science
Foundation Grant DMS--9401090 and wishes to thank the Mathematics Department
of the University of Melbourne for its hospitality while this paper 
was written.}

\author[Shapiro]{Michael Shapiro}

\thanks{The second author was partially supported by funds from the
Australian Research Council, the Group Theory Cooperative at City College
and the National Science Foundation; and he wishes to thank the Stekhlov
Institute for its hospitality.}

\begin{document}

\begin{abstract} Accessible groups whose word problems are
accepted by a \nsa\ are virtually free.  \end{abstract}

\maketitle

\section{Introduction.}\label{intro}

\noindent
During the past several years combinatorial group theory has received
an infusion of ideas both from topology and from the theory of formal
languages.  The resulting interplay between groups, the geometry of
their Cayley diagrams, and associated formal languages has led to
several developments including the introduction of automatic groups
\cite{E+}, hyperbolic groups \cite{Gr}, and geometric and
language--theoretic characterizations of virtually free groups
\cite{MS}. 

We will restrict our attention to finitely generated groups.  For any
such group the language of all words which define the identity is
called the word problem of the group.  By words we mean words over the
generators.  Of course the word problem depends on the choice of
generators.  In \cite{MS} virtually free groups are shown to be
exactly those groups whose word problem with respect to any set of
generators is a context--free language.  We are interested in
investigating groups whose word problems lie in other language
classes.

Formal languages are often defined in terms of the type of machine
which can tell whether or not a given word is in the language.  Such a
machine is said to accept the language.  Context--free languages are
accepted by pushdown automata, and those context--free languages which
are word problems are accepted by the subclass of deterministic
limited erasing pushdown automata~\cite[Lemma 3]{MS}. In this paper
we show that the more powerful class of deterministic limited erasing
nested stack automata accept exactly the same word problems.

\begin{theorem} \label{nsatheorem} Suppose $\G$ is an accessible group
whose word problem is recognized by a \nsa \ which accepts by final
state and empty stack; then $\G$ is virtually free.  \end{theorem}

The question of whether or not every group whose word problem is
accepted by a nested stack automaton is virtually free has been open
for some time.  Some possible counterexamples are proposed
in~\cite{LR}.  Notice our assumption that $\G$ is accessible.  While
it is not difficult to show that a group with context--free word
problem is finitely presented (and therefore accessible), deciding the
same question for groups whose word problem is accepted by a nested
stack automaton seems much harder.

In \cite{MS} and \cite{MS2} virtually free groups are characterized by
geometric conditions on their Cayley diagrams as well as by
language--theoretic conditions on their word problems.  In the course
of proving Theorem~\ref{nsatheorem} we are led to another such
geometric condition.  To express this condition we recall that a
choice of generators determines a word metric on a group. The distance
between two group elements $g$ and $h$ is the length of the shortest
word representing $g\inv h$. This metric is just the restriction to
vertices of the path metric in the corresponding Cayley diagram.
Different word metrics for the same group are quasi--isometric. (See
Definition~\ref{qidef}.) We will assume that every group is equipped
with a word metric. 

\begin{definition}\label{narrowdef} A group is \df{narrow} if there
exists an integer $i$ such that for any ball $B$ and all but finitely
many other balls $B'$ of the same radius, $B$ is separated from $B'$
by a set of size at most $i$.  A group which is not narrow is
\df{wide}.  \end{definition}

\noindent
Two subsets are \df{separated} by a set $S$ if every path between them
intersects $S$.  By \df{path} we mean a finite sequence of points each
a distance one from its successor. In other words a path is
just the sequence of vertices occurring in a path in the Cayley
diagram. 

Whether or not a group is narrow seems to depend on its word metric,
but in fact it does not.

\begin{theorem} \label{invarianttheorem} Narrowness is a quasi--isometry
invariant of groups. \end{theorem}

\noindent In adition narrowness characterizes those accessible groups 
which are virtually free.

\begin{theorem} \label{widenesstheorem} If the group $\G$ is 
finitely generated and accessible, then the following conditions are 
equivalent. 
\begin{enumerate} 
\item $\G$ is wide;
\item $\G$ is not virtually free;  
\item $\G$ contains a one--ended subgroup.  
\end{enumerate} 
If $\G$ is wide, then a one--ended subgroup can be found by splitting $\G$
over finite subgroups.  \end{theorem}

Narrowness is related to two geometric conditions shown in~\cite{MS2}
to characterize virtually free groups.  By~\cite[Theorem 2.9]{MS2} any
finitely generated group $\G$ (accessible or not) is virtually free if
and only if the components of the complements of all balls fall into
finitely many isomorphism classes of labeled graphs.  In fact by a
remark in the proof of that theorem, $G$ is virtually free if and only
if the frontiers of the components have uniformly bounded size.
Clearly this condition implies that $G$ is narrow, as any ball $B$ is
separated from all but finitely many balls of the same radius by the
frontiers of the components of the complement of $B$.  Conversely if
$G$ is narrow and accessible, then by~\cite{MS2} and
Theorem~\ref{widenesstheorem}, $G$ satisfies the other two geometric
conditions.  We do not know if there are narrow groups which are not
virtually free.

We are indebted to Swarup Gadde, Chuck Miller, Walter Neumann,
Christophe Pittet, and Nick Wormald for helpful conversations and in
particular to Pittet for an argument which appeared in earlier
versions of this paper.  The second author wishes to extend special
thanks to Rostislav Grigorchuk for his hospitality, and for first
introducing him to this problem.

The reader is referred to Aho's original paper~\cite{A2} on nested
stack automata, and for background on language theory to Hopcroft and
Ullman's book~\cite{HU} on automata theory.  The books by
Cohen~\cite{Co} and Lyndon and Schupp~\cite{LS} are references for
combinatorial group theory.

\section{Preliminaries.}

\noindent
\begin{definition} \label{qidef} A map $\phi:X\to X'$ between metric
spaces $(X,d)$ and $(X',d')$ is a \df{quasi--isometry} if there exists
a positive constant $k$ such that \begin{enumerate} \item $(1/k)d(x,y)
- k \le d(\phi(x),\phi(y)) \le kd(x,y) + k$; \item $X'=\cup_{x\in
X}B_k(\phi(x))$ \end{enumerate} \end{definition}

For example inclusion of a group with the word metric corresponding to
a choice of generators into the Cayley determined by the same
generators is a quasi--isometry. 

Quasi--isometry of metric spaces is an equivalence relation. As
different choices of finite generating set for a group yield
quasi--isometric metrics, it follows that all Cayley diagrams for a
finitely generated group are quasi--isometric. In fact Cayley diagrams
for commensurable groups are quasi--isometric~\cite[Prop.\ 11, page
8]{GH}.

Recall that all groups are assumed to be finitely generated. A
\df{virtually free} group is one with a free subgroup of finite index;
\df{virtually cyclic} groups are defined likewise.  All finite groups are 
virtually cyclic.

A group $\G$ is \df{one--ended} if for all $r$ the complement in its
Cayley diagram of the ball of radius $r$ around the identity has
exactly one infinite component. It is not hard to show that the
validity of this condition is independent of the generating set of the
group.  A famous theorem of Stallings~\cite{St} says that a group
which is not one--ended is either virtually cyclic, or splits as a free
product of two factors with a finite subgroup amalgamated, or is an HNN
extension with one stable letter and finite associated subgroups.

In the latter two cases we have $\G=H_1*_{K_0}H_2$ or $\G=\langle H_0,
t \mid t\inv K_1t=K_2\rangle$ respectively for subgroups $H_i$ of
$\G$.  We call $H_i$ a factor and $K_i$ an associated subgroup.  In
either case we say that $\G$ splits over a finite subgroup.

When a group $\G$ splits over a finite subgroup, it may be possible
that a factor splits over a finite subgroup, and one of the factors
of that splitting splits again etc. $\G$ is \df{accessible} if there
is an upper bound on the length of any such chain of splittings.  The
least upper bound is the accessibility length of $\G$.  It is a result
of Dunwoody~\cite{D} that a finitely presented group is accessible.

\begin{lemma} \label{one--ended} An accessible group is either 
virtually free or contains a one--ended subgroup. \end{lemma}

\begin{proof} If $\G$ does not split over a finite subgroup, then it is
virtually cyclic or one--ended.  If it does split, use induction on
accessibility length together with results of Gregorac~\cite{Gre} and
Karrass, Pietrowski and Solitar~\cite{KPS} which say that a group which
splits over a finite subgroup is virtually free if each factor is.
\end{proof}

\section{Machines.}

\noindent
We begin with an informal account based on the original definition of nested
stack automata as a certain kind of computer~\cite{A2}.  However this
definition is unwieldy, so our subsequent formal definition is in terms of
labeled graphs.  This approach to automata theory is well established.  See
for example Brainerd and Landweber~\cite[Chapter~4]{BL},
Eilenberg~\cite[Volume~A, Chapter~X]{Ei}, Floyd and Biegel~\cite{FB},
Gilman~\cite{Gi}, Goldstine~\cite{Go}, and Salomaa, Wood and Yu~\cite{SWY}.

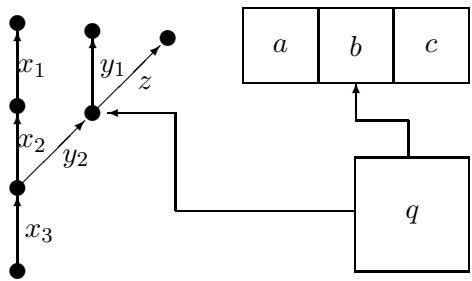
\begin{figure}[htb] \begin{center} 
\bigskip 
%\input{nsa.pic} % TexCAD file
%TexCad Options
%\grade{\off}
%\emlines{\off}
%\beziermacro{\off}
%\reduce{\on}
%\snapping{\on}
%\quality{2.00}
%\graddiff{0.01}
%\snapasp{1}
%\zoom{1.00}
\unitlength 1.00mm
\linethickness{0.4pt}
\begin{picture}(64.00,37.00)
\put(49.00,2.00){\framebox(15.00,15.00)[cc]{$q$}}
\put(56.00,17.00){\line(0,1){5.00}}
\put(56.00,22.00){\line(-1,0){7.00}}
\put(4.00,2.00){\circle*{2.00}}
\put(4.00,2.00){\vector(0,1){10.00}}
\put(4.00,13.00){\circle*{2.00}}
\put(4.00,13.00){\vector(0,1){10.00}}
\put(4.00,24.00){\circle*{2.00}}
\put(4.00,24.00){\vector(0,1){10.00}}
\put(4.00,13.00){\vector(1,1){9.00}}
\put(14.00,23.00){\circle*{2.00}}
\put(14.00,23.00){\vector(0,1){10.00}}
\put(14.00,23.00){\vector(1,1){9.00}}
\put(5.00,7.00){\makebox(0,0)[lc]{$x_3$}}
\put(4.00,19.00){\makebox(0,0)[lc]{$x_2$}}
\put(4.00,29.00){\makebox(0,0)[lc]{$x_1$}}
\put(15.00,29.00){\makebox(0,0)[lc]{$y_1$}}
\put(10.00,17.00){\makebox(0,0)[lc]{$y_2$}}
\put(20.00,27.00){\makebox(0,0)[lc]{$z$}}
\put(49.00,22.00){\vector(0,1){5.00}}
\put(4.00,35.00){\circle*{2.00}}
\put(14.00,34.00){\circle*{2.00}}
\put(24.00,33.00){\circle*{2.00}}
\put(25.00,10.00){\line(1,0){24.00}}
\put(34.00,27.00){\line(0,1){10.00}}
\put(34.00,37.00){\line(1,0){30.00}}
\put(64.00,37.00){\line(0,-1){10.00}}
\put(64.00,27.00){\line(-1,0){30.00}}
\put(44.00,27.00){\line(0,1){10.00}}
\put(54.00,27.00){\line(0,1){10.00}}
\put(39.00,32.00){\makebox(0,0)[cc]{$a$}}
\put(49.00,32.00){\makebox(0,0)[cc]{$b$}}
\put(59.00,32.00){\makebox(0,0)[cc]{$c$}}
\put(25.00,23.00){\line(0,-1){13.00}}
\put(25.00,23.00){\vector(-1,0){9.00}}
\end{picture}
\end{center} 
\caption{A nested stack automaton.\label{nsafigure}}
\end{figure}

Figure~\ref{nsafigure} shows a nested stack automaton $\A$ consisting of a
finite one--way input tape holding a word form an input alphabet $\Si$, a
finite set of internal states $\{q_i\}$, and a memory which at each instant
contains a finite directed tree with edge labels from a memory alphabet
$\Xi$. At each instant $\A$ is in a particular internal state, $q_i$, and is
either scanning a cell on the input tape or has moved off the tape to the
right.  Also $\A$ is pointing at a vertex of its memory tree. The label of 
the inedge to this vertex is called the current memory symbol. If $\A$ is 
pointing to the root, the current memory symbol is the empty word $\e$.

$\A$ begins a computation in a designated initial state $q_0$ and
scanning the leftmost cell on its input tape (or off the tape if the
input is $\e$).  Initially the memory is empty; that is, the memory
tree consists of just a single vertex.  Each computation consists of a
number of moves, and $\A$ is completely specified by a list of moves
to be made for various combinations of internal state, input letter,
and current memory symbol.  The input letter is either the content of
the current cell on the input tape or $\e$.  In the former case $\A$
moves right on the input tape as the last part of the move; in the
latter case it does not.  If $\A$ has moved off the input tape, only
$\e$--moves are possible. For any particular combination of internal state, 
input letter, and memory symbol there may be one, many or no moves 
specified.

If there is a sequence of moves in which $\A$ moves off the input tape and
reaches a configuration with empty memory and one of a designated set of
final states, then $\A$ \df{accepts} the word on the input tape.  The set of
all accepted words is the \df{language accepted by} $\A$.  If $\A$ reaches a
situation in which no move is possible, it halts.  If the conditions just
mentioned obtain, then the input is accepted, otherwise not.  $\A$ need not
halt after accepting an input, but its susequent behavior has no effect on 
the language accepted.

In addition to updating the pointer to the input tape, the other parts
of a move are the choice of a new internal state and a memory
operation.  There are four memory operations in addition to the
trivial operation $1$ in which the memory tree $\T$ is left unchanged.
To define these operations we observe that during a computation
vertices are added to and deleted from $\T$.  Thus at any given time
the vertices are ordered according to when they were added to $\T$.
The root is the earliest vertex and is never deleted.  The four memory
operations are moving the pointer down to the the parent of the
current vertex if that vertex is not the root, moving the pointer up to the 
latest child of a vertex if a child exists, deleting the
current vertex if it is a leaf but not the root, and adding a new edge whose
source is the current vertex and whose target is a new leaf.  After
deleting a leaf the current vertex is set to the source of the edge
to that leaf, and after adding a new edge the current vertex is set to
the new leaf.

It follows from the definitions above that the the vertex pointed to
by $\A$ is always on the path from the root to the latest vertex,
which is a leaf.  In particular only the latest vertex can be deleted.
Our description is not quite the same as the original in~\cite{A2}.
There the memory tree of Figure~\ref{nsafigure} would be replaced by a
memory tape containing three nested stacks $\$x_1x_2\$y_1\$z\cents
\hat{y_2}\cents x_3\#$.  In addition we allow NSA's to operate on
empty stack.  We leave it to the reader to check that the two kinds of
machine are equivalent in the sense that each can simulate the other.

Now we will give a more precise definition of nested stack automaton.
We fix once and for all an infinite countable memory alphabet $\Xi$.
Each NSA will use only finitely many letters from $\Xi$, so this
convention does no harm.  We order the vertices of a finite tree by
depth first search.  That is, the root is earliest followed by all the
vertices in order along a path to a leaf.  Then we back up to the
first vertex with an outedge which has not been traversed and continue
along that edge to another leaf, etc.  Edges are ordered according to
their target vertices.  If the tree in Figure~\ref{nsafigure} is
ordered this way by taking the leftmost possible outedge at each
opportunity, the corresponding ordering of edge labels would be
$x_3,x_2,x_1,y_2,y_1,z$.

It is clear that for any tree ordered by depth--first seach the latest
vertex is a leaf and deleting it gives a tree ordered in the same way.
Likewise adding a new edge with source anywhere along the path from
the root to the latest vertex and making the new leaf later than all
the other vertices in the tree yields a tree ordered in the right way.

\begin{definition}The set $\mtrees$ of \df{memory trees} consists of
all finite trees $\T$ ordered by depth--first search with 
\begin{enumerate} 
\item Root vertex $v_0$; 
\item Edges labeled by letters from $\Xi$; 
\item All edges directed away from the root; 
\item One distinguished vertex on the path from $v_0$ to the latest vertex
of 
$\T$. 
\end{enumerate}  
$\T_0$ is the tree consisting of just $v_0$.
\end{definition}

Next we define a monoid of operators on memory trees.

\begin{definition} \label{monoiddef} $\Mnsa$ is the monoid generated under 
composition by certain partial maps from $\mtrees$ to itself.  Pick $\T \in
\mtrees$ with distinguished vertex $v$, and let $y$ be the label of
the inedge to $v$.  If $v=v_0$, then $y=\e$.  We describe the effect
of the partial maps on $\T$.  In each case if $\T$ does not satisfy
the conditions given, then the map is not defined at $\T$.

\begin{description}

\item[$D_x(\T)$, $x\in \Xi$] If $x=y$ and $v \ne v_0$, then $D_x(\T)$ is 
obtained by changing the distinguished vertex of $\T$ to the parent of $v$.

\item[$U_x(\T)$, $x\in \Xi \cup\{\e\}$] If $x=y$ and $v$ is not a leaf, 
make the latest child of $v$ the new distinguished vertex.

\item[$P_x(\T)$, $x\in \Xi$] Add to $\T$ a new edge with source $v$, label
$x$, and target a new vertex $v_1$.  Make $v_1$ the latest vertex of $\T$
and the new current vertex.

\item[$Q_x(\T)$, $x\in \Xi$] If $x=y$ and $v$ is a leaf with parent 
$v_1$, delete $v$ and its inedge.  Make $v_1$ the distinguished vertex.

\end{description}
\end{definition}

\noindent Clearly $\Mnsa$ acts by partial injective maps whence it is
a submonoid of the symmetric inverse semigroup on $\mtrees$.  $\Mnsa$
has both an identity and a zero element, which we denote by $1$ and
$0$ respectively. 

\begin{definition} \label{nsa--def} Let $\Si$ be a finite alphabet.  A
\df{nested stack automaton} $\A$ over $\Si$ is a finite directed graph with
a designated initial vertex, designated final vertices, and with edges
labeled by pairs $(m,a)$ where $a\in \Si \cup \{\e\}$, and either $m=1$ or
$m$ is one of the generators defined in Definition~\ref{monoiddef}. 
\end{definition}

\noindent Every (directed) path in $\A$ has a \df{label} $(m,w)$
formed by multiplying the components of the edge labels in order.  A
path of length zero has label $(1,\e)$.

\begin{definition} \label{nsa--accept} A \df{computation} of a nested
stack automaton $\A$ is a path $\gamma$ which starts at the initial
vertex of $\A$ and has label $(m,w)$ for some $m$ with $\T=m(\T_0)$
defined.  $\T$ is called the \df{outcome} of $\gamma$.  The word
$w$ is \df{accepted} by $\A$ if there is a computation with label
$(w,m)$ and outcome $\T=\T_0$ ending at a final state. These computations 
are called \df{successful}. The set of all accepted words is the 
\df{language accepted} by $\A$. \end{definition}

\noindent In other words $\A$ accepts $w$ if it can read all of $w$,
empty its memory, and stop at a final state. One can also consider
NSA's which accept just by final state or empty memory, but we do not
do so here.

\begin{figure}[htb]
\bigskip
\centering
%\input{nsa3.pic}
%TexCad Options
%\grade{\on}
%\emlines{\off}
%\beziermacro{\on}
%\reduce{\on}
%\snapping{\on}
%\quality{2.00}
%\graddiff{0.01}
%\snapasp{1}
%\zoom{1.00}
\unitlength 0.80mm
\linethickness{0.4pt}
\begin{picture}(114.00,28.00)
\put(36.00,13.00){\circle{10.00}}
\put(36.00,17.00){\oval(6.00,20.00)[t]}
%\vector(39.00,21.00)(39.00,17.00)
\put(39.00,17.00){\vector(0,-1){0.2}}
\put(39.00,21.00){\line(0,-1){4.00}}
%\end
\put(36.00,28.00){\makebox(0,0)[cb]{$(P_x,a)$}}
\put(51.00,14.00){\makebox(0,0)[cb]{$(D_x,b)$}}
\put(36.00,13.00){\makebox(0,0)[cc]{$2$}}
\put(67.00,13.00){\circle{10.00}}
\put(98.00,13.00){\circle{10.00}}
\put(67.00,17.00){\oval(6.00,20.00)[t]}
\put(98.00,17.00){\oval(6.00,20.00)[t]}
%\vector(70.00,21.00)(70.00,17.00)
\put(70.00,17.00){\vector(0,-1){0.2}}
\put(70.00,21.00){\line(0,-1){4.00}}
%\end
%\vector(101.00,21.00)(101.00,17.00)
\put(101.00,17.00){\vector(0,-1){0.2}}
\put(101.00,21.00){\line(0,-1){4.00}}
%\end
\put(67.00,28.00){\makebox(0,0)[cb]{$(D_x,b)$}}
\put(98.00,28.00){\makebox(0,0)[cb]{$(U_x,c)$}}
\put(82.00,14.00){\makebox(0,0)[cb]{$(U_y,c)$}}
\put(67.00,13.00){\makebox(0,0)[cc]{$3$}}
\put(98.00,13.00){\makebox(0,0)[cc]{$4$}}
\put(5.00,13.00){\circle{10.00}}
\put(20.00,14.00){\makebox(0,0)[cb]{$(P_y,\e)$}}
\put(5.00,13.00){\makebox(0,0)[cc]{$1$}}
\put(102.00,13.00){\oval(20.00,6.00)[r]}
%\vector(106.00,10.00)(102.00,10.00)
\put(102.00,10.00){\vector(-1,0){0.2}}
\put(106.00,10.00){\line(-1,0){4.00}}
%\end
\put(114.00,13.00){\makebox(0,0)[lc]{$(Q_x,d)$}}
%\vector(10.00,13.00)(31.00,13.00)
\put(31.00,13.00){\vector(1,0){0.2}}
\put(10.00,13.00){\line(1,0){21.00}}
%\end
%\vector(41.00,13.00)(62.00,13.00)
\put(62.00,13.00){\vector(1,0){0.2}}
\put(41.00,13.00){\line(1,0){21.00}}
%\end
%\vector(72.00,13.00)(93.00,13.00)
\put(93.00,13.00){\vector(1,0){0.2}}
\put(72.00,13.00){\line(1,0){21.00}}
%\end
\put(51.00,1.00){\makebox(0,0)[cb]{$(Q_y,\e)$}}
%\emline(98.00,3.00)(98.00,8.00)
\put(98.00,3.00){\line(0,1){5.00}}
%\end
%\vector(5.00,3.00)(5.00,8.00)
\put(5.00,8.00){\vector(0,1){0.2}}
\put(5.00,3.00){\line(0,1){5.00}}
%\end
\put(51.50,3.00){\oval(93.00,6.00)[b]}
\end{picture}
\caption{An NSA $\A$ which accepts $\{a^nb^nc^nd^n\}^*$. \label{nsagraph}}
\end{figure}
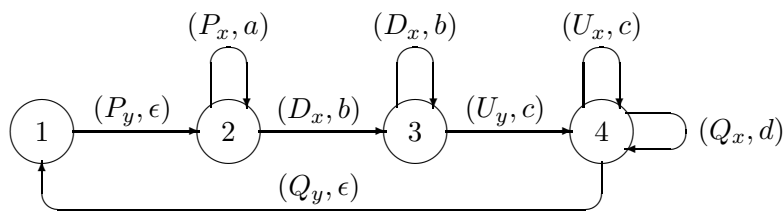

Figure~\ref{nsagraph} shows an NSA $\A$ which accepts $\{a^nb^nc^nd^n\}^* =
\{(a^nb^nc^nd^n)^k \mid n\ge 0, k\ge 0\}$.  Vertex 1 is both the start
vertex and the single final vertex.  To see that $\A$ accepts the language
claimed, first observe that domain of $U_x$ is disjoint from the range of
$Q_x$.  Thus a computation in the sense of Definition~\ref{nsa--accept}
cannot have an edge with label $(Q_x,d)$ followed by one with label
$(U_x,c)$.  Consequently the label of any successful computation by $\A$ is
a product of terms of the form $(P_yP_x^iD_x^jU_yU_x^kQ_x^lQ_y,
a^ib^jc^kd^l)$.  In particular there is a successful computation whose
label is the empty product, $(1,\e)$, so $\A$ accepts $\e$.

To show that the language accepted by $\A$ is as claimed, it suffices to
show that the exponents $i,j,k,l$ in any term above are equal and that all
cases in which the exponents are equal and greater than $0$ can occur.
Suppose $(P_yP_x^iD_x^jU_yU_x^kQ_x^lQ_y, a^ib^jc^kd^l)$ is the first term.
$\A$ begins by constructing a memory tree $\T$ with a single branch labeled
$yx^i$.  If $i < j$, then $\A$ tries unsuccessfully to move down to the root
of $\T$ by executing $D_x$ while pointing at the vertex of $\T$ with inedge
labeled $y$.  If $i > j$, then $\A$ tries to move up past the leaf of $\T$.
Conseqently $\A$ reaches vertex $4$ if and only if $i=j>0 $ and $k>0$.
Additional arguments of a similar nature demonstrate that $\A$ returns to
vertex $1$ if and only if $i=j=k=l>1$ and that upon its return $\T=T_0$. 
But $\T=T_0$ imples that our analysis applies to each term in succession.

It follows easily from the preceding paragraph that at any point in a
computation the memory tree $\T$ has only a single branch.  In other words
$\A$ is a \df{stack automaton}; operations $P_x$ are executed only when the
distinguished vertex is a leaf.  In addition it is clear from
Figure~\ref{nsagraph} that $\A$ has at most one move for each combination of
state and input symbol.  Such automata are called deterministic.

\begin{definition} \label{deterministic} An NSA $\A$ over $\Si$ is
\df{deterministic} if each combination of vertex $v$, memory $\T$, and
letter $a\in \Si$ admits at most one outedge with label $(m,a)$ or
$(m,\e)$ such that $m(\T)$ is defined. In other words the $m$'s which
occur in labels of these outedges have pairwise disjoint domains in
$\mtrees$.  \end{definition}

\noindent Any computation of length $n$ by a deterministic NSA can be
continued in at most one way to a computation of length $n+1$.

\begin{definition} \label{limited} An NSA $\A$ has \df{limited erasing} if
there is a constant $k$ such that every path in $\A$ with label  $(m, \e)$ 
 has at most $k$ edges with labels involving $Q$. 
In other words the $\A$ makes at most $k$ erasures from its memory without 
consuming input. \end{definition}

\noindent  The NSA in Figure~\ref{nsagraph} is deterministic and has limited 
erasing with $k=1$.

\begin{lemma}\label{inverseh} Suppose $L\subseteq \SSi$ is the language of
all words accepted by the deterministic NSA with limited erasing $\A$ over
$\Si$.  If $\De$ is a finite alphabet and $f:\DD \to \SSi$ is a homomorphism
which does not map any generator to the empy word, then $f\inv(L)$ is also
accepted by a deterministic NSA with limited erasing. \end{lemma}
 
\begin{proof} Suppose first that $f$ maps generators to generators.
In this case we can construct the required NSA $\A' $ directly from
$\A$.  Replace each edge with label $(m,a)$, $a\in \Si$, by a set of
edges with labels $(m,p)$ for all $p\in f\inv(a)$.  The new edges have
the same source and target as the edge they replace.  If $f\inv(a)$ is
empty, the original ledge is simply deleted.  It is clear from the
construction that $\A' $ is deterministic with limited erasing.
Further for any path in $\A$ from vertex $v_1$ to $v_2$ with label
$(m,w)$ there is for each $w'\in f\inv(w)$ a path in $\A' $ from $v_1$
to $v_2$ label $(m,p)$.  Conversely every path in $\A' $ with label
$(m,w')$ projects to a path in $\A$ with label $(m,f(w'))$.  It
follows that $\A' $ accepts $f\inv(L)$.

In general $f$ factors as a product of homomorphisms of the type just
considered and homomorphisms which map one generator into a word of length
two and all other generators into themselves.  Thus it suffices to prove
$f\inv(L)$ is an NSA language in this second case.  More precisely we may
assume that $\De$ and $\Si$ are the same except that one generator $a \in
\De$ is replaced by $a_1$ and $a_2$ in $\Si$.  Further $f(a)=a_1a_2$ and 
$f$ maps all other elements of $\De$ to themselves.

Construct an NSA $\A' $ from the union of two disjoint copies of $\A$,
say $\A_1$ and $\A_2$. The idea here is to modify $\A_1\cup \A_2$ so that
paths in $\A$ correspond to paths in $\A_1$ except that subpaths with
label sequences 

$$(m_1,a_1), (m_2,\e), \ldots, (m_{k-1},\e), (m_k,a_2)$$

\noindent
move over to $\A_2$ at the first edge and return to $\A_1$ at the last
edge. Further the first label of such a subpath is changed to
$(m_1,a)$ and the last to $(m_k,\e)$. It is also necessary to prevent
$\A' $ from emptying its memory at any vertex in $\A_2$. For this
purpose we add a new memory symbol $z$ which is pushed on to the
memory at the beginning of each computation and which can only be
removed at final vertices in $\A_1$.

More precisely $\A' $ is formed from the disjoint union $\A_1 \cup
\A_2 \cup \{v_0,v_1\}$ in the following way. 
\begin{enumerate}
\item The start vertex of $\A' $ is $v_0$, and there is an edge from $v_0$ 
to the start vertex of $\A_1$ with label $(P_z,\e)$;
\item Every final vertex of $\A_1$ has an outedge to $v_1$ with label
$(Q_z,\e)$;
\item The single final vertex is $v_1$.
\item Each edge in $\A_1$ with $a_1$ in its label has $a_1$ changed to
$a$ and its terminal vertex changed to the corresponding vertex in $N_2$; 
\item All edges from $\A_2$ except those with an $\e$ or $a_2$ in their 
label are removed;
\item Each edge in $\A_2$ whose label involves $a_2$ has $a_2$ changed 
to $\e$ and  its terminal vertex changed to the corresponding vertex of 
$N_1$.
\end{enumerate}
It is straightforward to check that $\A' $ is deterministic with limited 
erasing, and accepts $f\inv(L)$.  \end{proof}

\section{Group languages.}

\noindent
In this section we develop the first properties of groups whose word
problems are solvable by deterministic nested stack automata with limited
erasing.  From now on NSA will refer to an anutomaton of this type which
accepts by final state and empty stack, and NSA language will mean a
language accepted by an NSA.

We fix some notation.  A {{\em choice of generators}} for a group $\G$
is a surjective monoid homomorphism $\si:\SSi\to G$ from a finitely
generated free monoid.  We will write $\ovr w$ for $\si(w)$ and assume
that a choice of generators $\si:\SSi\to \G$ always has {{\em formal
inverses}}.  That is, $\Si$ is a union of pairs $\{a,a\inv\}$ and
$\si(a\inv)=(\si(a))\inv$.  We emphasize that $\Si$ still generates
$\SSi$ freely as a monoid; there is no cancellation in $\SSi$.  Recall
that $\e$ denotes the empty word and $\ovr{\e}=1$. The {{\em word problem}} 
of $\G$ corresponding to a certain choice of generators is $\{w\in\SSi \mid
\ovr w = 1\}$.

\begin{lemma} If the word problem of a group $\G$ with respect to one choice
of generators is an NSA language, then so are the word problem with respect
to any choice of generators and the word problem for every finitely
generated subgroup of $\G$.  \end{lemma}

\begin{proof} Let $\si:\SSi\to \G$ be a generating set for $G$ such
that an NSA solves the word problem for $G$ with respect to $\SSi$,
and $\de:\DD\to H\subseteq \G$ be a generating set for $H$.  Choose a
homomorphism $f:\SSi\to\DD$ such that $\de\circ f = \si$ and $f$ does
not map any generator to the identity.  Apply Lemma~\ref{inverseh}.
\end{proof}

Now let $\A$ be an NSA over $\Si$ accepting the word problem of $\G$
with respect to a choice of generators $\SSi\to G$.  We use $\A$ to
construct a graph $\CG$ which covers both $\A$ and $\CD$, the Cayley
diagram of $\G$. We augment $\CD$ by adding an edge with label $\e$
from every vertex to itself. 

\begin{definition} A \df{configuration} of an NSA $\A$ is a pair $(q,\T)$
where $q$ is a state of $\A$ and $\T$ is a memory tree.  A configuration is
\df{accessible} if $\T$ is the outcome of a valid computation $\gamma$
ending at $q$ and if there is a continuation $\gamma'$ such that
$\gamma\gamma'$ is successful.  \end{definition}

Accessible configurations might more properly be called accessible and
co--accessible.

\begin{definition} The \df{configuration graph} $\CG$ of $\A$ has as
vertices all accessible configurations.  There is an edge from
$(q,\T)$ to $(q',\T')$ with label $a\in \Si\cup \{\e\}$ if and only if
there is an edge with label $(m,a)$ from $q$ to $q'$ in $\A$ such that
$m(T)=\T' $. The initial vertex of $\CG$ is $(q_0,T_0)$ where $q_0$ is
the initial vertex of $\A$.  \end{definition}

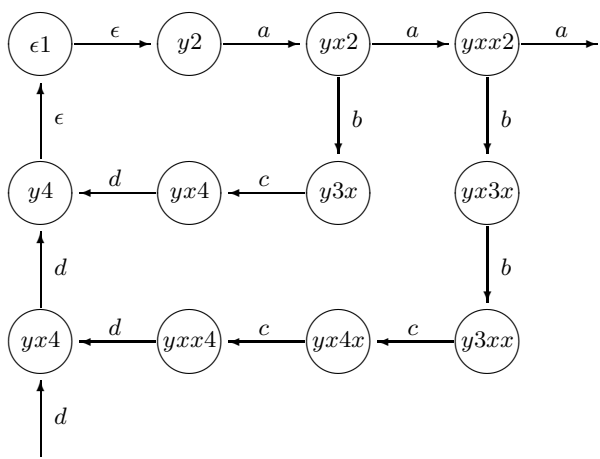
\begin{figure}[htb] \begin{center} 
\bigskip 
\footnotesize
%\input{nsa4.pic} % TexCAD file
%TexCad Options
%\grade{\on}
%\emlines{\off}
%\beziermacro{\on}
%\reduce{\on}
%\snapping{\on}
%\quality{2.00}
%\graddiff{0.01}
%\snapasp{1}
%\zoom{1.00}
\unitlength 0.90mm
\linethickness{0.4pt}
\begin{picture}(90.00,68.00)
\put(8.00,19.00){\circle{10.00}}
\put(8.00,41.00){\circle{10.00}}
\put(8.00,19.00){\makebox(0,0)[cc]{$yx4$}}
\put(8.00,41.00){\makebox(0,0)[cc]{$y4$}}
%\vector(25.00,19.00)(14.00,19.00)
\put(14.00,19.00){\vector(-1,0){0.2}}
\put(25.00,19.00){\line(-1,0){11.00}}
%\end
%\vector(25.00,41.00)(14.00,41.00)
\put(14.00,41.00){\vector(-1,0){0.2}}
\put(25.00,41.00){\line(-1,0){11.00}}
%\end
\put(30.00,19.00){\circle{10.00}}
\put(30.00,41.00){\circle{10.00}}
\put(30.00,63.00){\circle{10.00}}
\put(52.00,19.00){\circle{10.00}}
\put(52.00,41.00){\circle{10.00}}
\put(52.00,63.00){\circle{10.00}}
\put(74.00,19.00){\circle{10.00}}
\put(74.00,41.00){\circle{10.00}}
\put(74.00,63.00){\circle{10.00}}
\put(30.00,19.00){\makebox(0,0)[cc]{$yxx4$}}
\put(30.00,41.00){\makebox(0,0)[cc]{$yx4$}}
\put(30.00,63.00){\makebox(0,0)[cc]{$y2$}}
\put(52.00,19.00){\makebox(0,0)[cc]{$yx4x$}}
\put(52.00,41.00){\makebox(0,0)[cc]{$y3x$}}
\put(52.00,63.00){\makebox(0,0)[cc]{$yx2$}}
\put(74.00,19.00){\makebox(0,0)[cc]{$y3xx$}}
\put(74.00,41.00){\makebox(0,0)[cc]{$yx3x$}}
\put(74.00,63.00){\makebox(0,0)[cc]{$yxx2$}}
%\vector(47.00,19.00)(36.00,19.00)
\put(36.00,19.00){\vector(-1,0){0.2}}
\put(47.00,19.00){\line(-1,0){11.00}}
%\end
%\vector(69.00,19.00)(58.00,19.00)
\put(58.00,19.00){\vector(-1,0){0.2}}
\put(69.00,19.00){\line(-1,0){11.00}}
%\end
%\vector(47.00,41.00)(36.00,41.00)
\put(36.00,41.00){\vector(-1,0){0.2}}
\put(47.00,41.00){\line(-1,0){11.00}}
%\end
%\vector(35.00,63.00)(46.00,63.00)
\put(46.00,63.00){\vector(1,0){0.2}}
\put(35.00,63.00){\line(1,0){11.00}}
%\end
%\vector(57.00,63.00)(68.00,63.00)
\put(68.00,63.00){\vector(1,0){0.2}}
\put(57.00,63.00){\line(1,0){11.00}}
%\end
%\vector(79.00,63.00)(90.00,63.00)
\put(90.00,63.00){\vector(1,0){0.2}}
\put(79.00,63.00){\line(1,0){11.00}}
%\end
\put(41.00,64.00){\makebox(0,0)[cb]{$a$}}
\put(63.00,64.00){\makebox(0,0)[cb]{$a$}}
\put(85.00,64.00){\makebox(0,0)[cb]{$a$}}
\put(41.00,42.00){\makebox(0,0)[cb]{$c$}}
\put(41.00,20.00){\makebox(0,0)[cb]{$c$}}
\put(63.00,20.00){\makebox(0,0)[cb]{$c$}}
\put(19.00,42.00){\makebox(0,0)[cb]{$d$}}
\put(19.00,20.00){\makebox(0,0)[cb]{$d$}}
\put(10.00,8.00){\makebox(0,0)[lc]{$d$}}
\put(10.00,30.00){\makebox(0,0)[lc]{$d$}}
\put(76.00,30.00){\makebox(0,0)[lc]{$b$}}
\put(76.00,52.00){\makebox(0,0)[lc]{$b$}}
\put(54.00,52.00){\makebox(0,0)[lc]{$b$}}
%\vector(8.00,2.00)(8.00,13.00)
\put(8.00,13.00){\vector(0,1){0.2}}
\put(8.00,2.00){\line(0,1){11.00}}
%\end
%\vector(8.00,24.00)(8.00,35.00)
\put(8.00,35.00){\vector(0,1){0.2}}
\put(8.00,24.00){\line(0,1){11.00}}
%\end
%\vector(52.00,58.00)(52.00,47.00)
\put(52.00,47.00){\vector(0,-1){0.2}}
\put(52.00,58.00){\line(0,-1){11.00}}
%\end
%\vector(74.00,58.00)(74.00,47.00)
\put(74.00,47.00){\vector(0,-1){0.2}}
\put(74.00,58.00){\line(0,-1){11.00}}
%\end
%\vector(74.00,36.00)(74.00,25.00)
\put(74.00,25.00){\vector(0,-1){0.2}}
\put(74.00,36.00){\line(0,-1){11.00}}
%\end
\put(8.00,63.00){\circle{10.00}}
\put(8.00,63.00){\makebox(0,0)[cc]{$\e1$}}
\put(10.00,52.00){\makebox(0,0)[lc]{$\e$}}
%\vector(8.00,46.00)(8.00,57.00)
\put(8.00,57.00){\vector(0,1){0.2}}
\put(8.00,46.00){\line(0,1){11.00}}
%\end
%\vector(13.00,63.00)(24.00,63.00)
\put(24.00,63.00){\vector(1,0){0.2}}
\put(13.00,63.00){\line(1,0){11.00}}
%\end
\put(19.00,64.00){\makebox(0,0)[cb]{$\e$}}
\end{picture}
\normalsize
\end{center} 
\caption{Part of the configuration graph of the NSA $\A$ of
Figure~\ref{nsagraph}.\label{configurationgraph}} \end{figure}

In Figure~\ref{configurationgraph} a vertex with label $yx3x$, say, stands
for the configuration $(3,\T)$ in which $\T$ consists of one branch of
length three with label $yxx$ and distinguished vertex a distance two from
the root.

The the first half of the next lemma is clear from 
Definition~\ref{nsa--def} and Definition~\ref{nsa--accept}; the second half 
follows from Definitions~\ref{deterministic} and~\ref{limited}.

\begin{lemma}\label{configurationlemma} The following conditions hold.
\begin{enumerate}
\item Every computation $\gamma$ of $\A$ which can be continued to a 
successful computation lifts uniquely to a
path in $\CG$ which starts at $(q_0,\T_0)$ and has label equal to the 
first component of the label of $\gamma$.
\item If a computation which ends at $q$ with outcome $\T$ lifts to a path
starting at $(q_0,\T_0)$, then the lift ends at $(q,\T)$.  Conversely any 
path in $\A$ from $(q_0,\T_0)$ to a vertex $(q, \T)$ is the lift of a
computation with outcome $\T$.  
\item Each vertex of $\CG$ either has a single outedge with label $\e$ and 
no other outedges, or it has no outedges with label $\e$ and at most one 
outedge with label $a$ for each $a\in \Si$.
\item For some constant $K$ any
path in $\CG$ has at most $K$ successive edges labeled $\e$.
\end{enumerate} \end{lemma}

\begin{lemma} \label{covering} 
There is a homomorphism of labeled graphs $\phi:\CG \to \CD$, where
$\CD$ is Cayley diagram of $\G$ augmented by the addition of a loop
with label $\e$ at every vertex.  The image of a vertex in $\CG$ is
the group element represented by the label of any path from
$(q_0,T_0)$ to that vertex.  In particular the initial vertex of $\CG$
maps to $1$.  Every path $\CD$ starting at $1$ lifts to a path in
$\CG$ which starts at $(q_0,\T_0)$ and projects to a path differing
from the original only by addition or deletion of edges with label
$\e$.  The lift is unique up to a terminal segment with label $\e$.
\end{lemma}

\begin{proof} By definition of $\CG$ each vertex $(q,\T)$ is reached
by a path from $(q_0,\T_0)$.  Suppose there are two such paths with
labels $w$ and $w'$.  By the definition of $\CG$ again there is a path
from $(q,\T)$ to some vertex $(q',\T_0)$.  Let $u$ be the label of
this path.  Then $\A$ accepts $wu$ and $w'u$ whence both denote the
identity in $\G$.  It follows that $w$ and $w'$ represent the same
element of $\G$.

Define $\phi$ by mapping each vertex of $\CG$ to the group element
represented by the label of any path from $(q_0,\T_0)$ to that vertex.
In particular $\phi((q_0,\T_0))=1$.  Suppose there is an edge from
$(q,\T)$ to $(q',\T')$ with label $a\in \Si\cup\{\ \epsilon\}$.  Pick
a path from $(q_0,\T_0)$ to $(q,\T)$, and let its label be $w$.  As
$wa$ is the label of a path to $(q',\T')$, it follows that if
$\phi((q,\T))=g$, then $\phi((q',\T'))=g\ovr a$.  Thus $\phi$ is a
graph homomorphism.

Since $\phi$ is a homomorphism, the penultmate assertion amounts to showing 
that for every $w\in \SSi$ there is a path in $\CG$ starting at $(q_0,\T_0)$
and with label $w$. But for some $v\in \SSi$, $wv$ represents the identity 
in $\G$ and hence there is a successful computation in $\A$ with label 
$(m,w)$. The last assertion follows from the third part of 
Lemma~\ref{configurationlemma}. \end{proof}

We pause to remark on a difference in configuration graphs of NSA's
and pushdown automata which illustrates the different power of the two
types of automaton.  A pushdown automaton is a nested stack automaton,
not necessarily quasi-realtime or deterministic, which has no labels
involving $D_x$ or $U_x$.  In other words its memory trees all have
just one branch and the distinguished vertex is always the leaf.  We
metrize graphs in the usual way by disregarding orientation and taking
edges to be isometric to the unit interval.

\begin{theorem}\label{treequotient}
The configuration graph of a pushdown automaton is quasi-isometric to a
tree.  This quasi-isometry is given by a quotient map.  In contrast to this,
the configuration graph of a stack automaton (and hence an NSA) may have
arbitrarily large isometrically embedded loops.  \end{theorem}

\begin{proof} The second claim is clear from
Figure~\ref{configurationgraph}.  It remains to prove that the
configuration graph of a pushdown automaton is quasi--isometric to a
tree.  Let $\CG$ be the configuration graph of a pushdown automaton
$\A$.  We will write $(p,\T) \sim (q,\T)$ if there is an undirected
path in $\CG$ between $(p,\T)$ and $(q,\T)$ with the property that
$\T$ is an initial segment of the memory tree of every intermediate
vertex.  That is to say, $(p, \T) \sim (q, \T)$ if we can get from
$(p, \T)$ to $(q, \T)$ by a sequence of forward and backwards moves
without ever erasing any portion of $\T$.  

Since $\CG$ has no labels involving $D_x$, the existence of such a path
depends only on $p$, $q$, and the label of the inedge to the leaf of $\T$.
It follows that there is a universal bound on the undirected distance
between $(p, \T)$ to $(q,\T)$ whenever $(p,\T) \sim (q,\T)$.  We denote the
equivalence class of $(p,\T)$ by $([p],\T)$ since $\T$ is constant
throughout the class.  This is an abuse of notation, since $[p]$ depends on
$\T$.  We console ourselves with the hope that $\T$ will be clear from
context.

Project $\CG$ to the graph $\overline{\CG}$ whose vertices are
equivalence classes of $\sim$ with distinct classes joined by an
unoriented edge if there is an edge between any two vertices in the
preimages of the classes.  In the case where an edge of $\CG$ connects
a vertex to itself, we project that edge of $\CG$ to the image of its
endpoints in $\ovr{\CG}$.  We have seen that there is a universal
bound on the size of $([p], \T)$.  It follows that the quotient map
is a quasi--isometry.  We will show that $\ovr{\CG}$ has no simple
cycles and is therefore a tree.  

Suppose to the contrary that there is a simple cycle in $\ovr\CG$.  Since
$\ovr\CG$ has no loops and at most one edge between any two vertices, this
cycle must have length at least three.  We examine the cycle at a place
where $\T$ has maximal size.  Either there is a single edge between
distinct vertices $([p],\T)$ and $([q],\T)$, or there are two edges
connecting distinct vertices $([p],\T')$, $([q],\T')$ with $\T'$ a proper
prefix of $T$ to a vertex $([r],\T')$. The vertices $([p],\T')$ and  
$([q],\T')$ have the same memory tree because their memory trees are shorter
than $\T$ but derived from $\T$ by popping a single element.

In the first case, there is an edge in $\CG$ between $(p',\T)$ and
$(q',\T)$ for some $p'\in [p]$ and $q' \in [q]$.  Consequently
$([p],\T) = ([q], \T)$, and we have a contradiction. In the second
case there is an undirected path from $(p',\T')$ to $(q',\T')$ along
which each memory tree is an initial segment of $\T$.  Hence
$([p],\T') = ([q],\T')$, and our cycle is not simple.  It follows that
$\ovr\CG$ is a tree.  \end{proof}

Notice that the proof above does not work for the usual kind of pushdown
automata which allowed to push a string of letters onto the stack on one
move.  Our pushdown automata are restricted to pushing one symbol at a time.
However is easy to simulate a more general automaton by adding internal
states to obtain one of ours, and it is straightforward to show that the
configuration graphs of the two machines are quasi--isometric.  Thus the
result above holds for standard pushdown automata.

Theorem \ref{treequotient} provides another way of proving the result
of Muller and Schupp.  For we now have the configuration graph $\CG$
mapping to the Cayley graph $\CD$ and quasi--isometrically to a tree
$T$.  The lifting properties quickly show that $\CD$ is not one--ended
and one proceeds as before.  This points up a difference between our
result and that of Muller and Schupp.  In the case of pushdown
automata, the treeness --- and therefore, the freeness --- are already
implicit in the class of machines.  In the case of nested stack
automata, the freeness is a result of the class of automata together
with the fact that they are being used to solve the word problem.

\section{Proof of Theorem~\ref{invarianttheorem} }

\noindent
Suppose that $\phi:G \to \G'$ is a quasi--isometry of groups and that
$\G'$ is narrow.  It suffices to show that $\G$ is narrow too. Pick a
ball $B \subseteq \G$ of radius $r$ and let $\overline B$ be its image
in $\G'$.  By Definition~\ref{qidef} $\overline B$ lies in a ball $C$
of radius $k(r+1)$. Since $\G'$ is narrow, there is an integer $i$
such that all but finitely many balls $C'_j\subseteq \G'$ of radius
$k(r+1)$ are separated from $C$ by $i$ points.

It follows from Definition~\ref{qidef} that $\phi$ is uniformly finite
to one.  Consequently the preimage in $\G$ of the union of the balls
$C'_j$ not satisfying the separation condition is finite.  Thus with a
finite number of exceptions every ball $B'\subseteq \G$ of radius $r$
has an image which lies in a ball $C'$ separated from $C$ by a set $S$
of size $i$.

The image under $\phi$ of any path from $B$ to $B'$ is a sequence of
points from $C$ to $C$ with each point a distance at most $2k$ from
its successor.  Add points to obtain a path from $C$ to $C'$ and
observe that this path must intersect $S$. It follows that at least
one of the original image points is a distance at most $k$ from $S$.
In other words the preimage in $\G$ of all points within $k$ of $S$
disconnects $B$ and $B'$. Since $\phi$ is uniformly finite to one, we
are done.

\section{Proof of Theorem~\ref{widenesstheorem}}

\noindent
We must prove the following three implications:
\begin{enumerate}
\item  A group which is wide is not virtually free.
\item  A group which is not virtually free has a one--ended subgroup.
\item  A group with a one--ended subgroup is wide.
\end{enumerate}

The first of these follows from the discussion immediately following
the statement of the Theorem, and the second is a consequence of
Lemma~\ref{one--ended}.  Thus it remains only to prove that if $\G$
contains a one--ended subgroup, then it is wide. Choose generators for
$\G$ which contain generators of the subgroup; then it follows
immediately from Definition~\ref{narrowdef} that if $\G$ is narrow, so
is the subgroup. Conversely if the subgroup is wide, so is $\G$. Thus
it suffices to prove that a one--ended group is wide.

Suppose to the contrary that $\G$ is one--ended and narrow.  Assume we can
choose the sets $S$ in Definition~\ref{narrowdef} so that the distance
between any two points in $S$ is uniformly bounded.  It follows from the
fact that $\G$ acts transitively and isometrically on itself by left
translation that up to the action of $\G$ there are only finitely many
different $S$'s.  Thus at least one $S$ separates pairs of balls of
arbitrarily large radius whence $\G-S$ has at least two infinite connected
components contrary to the hypothesis that $\G$ is one--ended.

We complete the proof by showing that we can choose the sets $S$ as
required.  Observe that since $\G$ is one--ended and
of bounded valence, $\G-\{1\}$ has one infinite component and finitely many
finite ones.  Some ball $C$ of radius at least 1 around $1$ contains all the 
finite components, and any two vertices on the boundary of $C$ are joined by
a path lying entirely in the infinite component. Pick one such path
for each pair of vertices on the boundary of $C$, and pick a ball $D$
around $1$ containing $C$ and all these paths. Let $D$ have radius
$d$. It is clear from Definition~\ref{narrowdef} that we may 
assume the balls to be separated are of radius greater than $d$.

It does no harm to restrict attention to sets $S$ which are minimal
with respect to inclusion, and we do so.  Pick balls $B$ and $B'$ of
radius greater than $d$ and separated by a set $S$.  Suppose $S$
contains a point $v$ at a distance greater than $d$ from all other
points of $S$.  Without loss of generality we assume this vertex is
$1$ whence $D\cap S = \{1\}$.  

By minimality of $S$ there is a path $\gamma$ from from $B$ to $B'$ with
$\gamma \cap S=\{1\}$.  If $\gamma$ starts outside $C$, define $\gamma_1$ to
be the initial segment of $\gamma$ from $B$ to a point $x$ on the boundary
of $C$.  Otherwise $B\cap C$ is not empty; and because $B$ is connected and
larger than $C$, there is a point $x$ in $B$ and on the boundary of $C$.  In
this case define $\gamma_1$ to be the path of length zero from $x$ to
itself.  In both cases $\gamma_1$ goes from $B$ to a point $x$ on the
boundary of $C$ and intersects $S$ trivially.  Define $\gamma_3$ similarly
from a point $y$ on the boundary of $C$ to $B'$, and let $\gamma_2$ be a
path from $x$ to $y$ lying in $D-\{1\}$.  Clearly $\gamma_1\gamma_2\gamma_3$
is a path from $B$ to $B'$ intersecting $S$ trivially.  Hence $S$ cannot
have a point $v$ a distance greater than $d$ from all other points of $S$. 
Since the sets $S$ have uniformly bounded size, we are done.

\section{Proof of Theorem~\ref{nsatheorem}}

\noindent
Suppose that the word problem of $\G$ is accepted by a \nsa\ $\A$ with
$q$ states and limited erasing constant $K$.  To prove
Theorem~\ref{nsatheorem} it is enough to show that $\G$ is narrow.
Recall that the configuration graph $\CG$ of $\A$ projects onto the
Cayley diagram $\CD$ of $\G$ where $\CD$ is augmented by the addition
of edge with label $\e$ from every vertex to itself.  For any
constant $C$ there are only finitely many vertices in $\CG$ with
memory tree of at most $C$ edges.  Thus for any element $g\in G$ far
enough away from $1$ every vertex $\CG$ projecting to $g$ has memory
tree with at least $C$ edges.

By the action of $\G$ we may take $B$ to be a ball around $1$; and by
discarding only finitely many possiblities for $B'$ we may assume that
the center of $B'$ is an element $g$ far enough away from $1$ so that
every vertex in the preimage in $\A$ of $g$ has a memory tree with at
least $2Kr+1$ edges.  If necessary move $g$ farther away so that $B$
and $B'$ are a distance at least 2 from each other.

Fix a path $\gamma$ from $1$ to $g$ in $\CD$.  By Lemma~\ref{covering} there
is a unique shortest path $\hat{\gamma}$ beginning at the initial vertex
$(q_0,\T_0)$ of $\A$ and projecting to $\gamma$.  By choice of $g$,
$\hat{\gamma}$ ends at $(q,\T)$ for some memory tree with at least $2Kr+3$
edges.  Recall that the edges of $\T$ are ordered, and let $\T'$ be the tree
obtained by removing the $Kr+1$ latest edges and taking the latest remaining
vertex as distinguished.  We claim that the images in $\G$ of all vertices
of the form $(q',\T')$ separate $B$ and $B'$.

It suffices to show that any path beginning at the boundary of $B'$
and ending at the boundary of $B$ contains one of the desired
vertices. Extend such a path by geodesics of length $r$ to a path
$\gamma_1$ from $g$ to $1$ in $G$, and consider the cycle
$\gamma\gamma_1$.  As the label $w$ of this cycle represents the
identity in $\G$, there must be a path in $\CG$ from $(q_0,\T_0)$ to a
vertex $(q_1,\T_0)$ with label $w$ and $q_1$ a final state.  Our choice
of $\hat{\gamma}$ insures that $\hat{\gamma}$ is a prefix of this path.
Consequently there is a path $\hat{\gamma_1}$ from $(q,\T)$ to $(q_1,\T_0)$
projecting to $\gamma_1$.

Since $\A$ is limited erasing, the initial and terminal segments of
$\hat{\gamma_1}$ projecting to the geodesic segments at the ends of
$\gamma$ have length at most $Kr$. It follows at most $Kr+1$ edges of
$\T$ can be deleted along each of these initial and terminal segments. Since
edges must be deleted in order, it follows that there is an edge not
in these segments with label $(q',\T')$.

\bibliographystyle{amsplain}

\begin{thebibliography}{99}

\bibitem{A2} A. Aho, Nested Stack Automata,
J. Assoc.\  Computing Machinery, {\bf 16} 1969, 383--406.

\bibitem{BL}
W. Brainerd and L. Landweber, {\em Theory of Computation},
John Wiley \& Sons, New York, 1974.

\bibitem{Co} D. E. Cohen, Combinatorial Group Theory: A Topological 
Approach, London Mat.\ Soc.\ 1989.

\bibitem{D} M. Dunwoody, The accessibility of finitely presented groups,
Invent.  Math.  {\bf 81} 1985, 449--457.

\bibitem{Ei}
S. Eilenberg, {\em Automata, Languages and Machines}, vols. A and B, 
Academic Press, New York, 1974.

\bibitem{E+} D.B.A.  Epstein, J.W.  Cannon, D.F.  Holt, S.V.F.  Levy, M.S.
Paterson, W.P.  Thurston, {\em Word processing in groups}, Jones and
Bartlett,
1992

\bibitem{FB}
R. Floyd and R. Biegel, {\em The Language of Machines}, Computer Science 
Press, New York, 1994.

\bibitem{Gi} R. Gilman, Formal languages and infinite groups, in {\em
Geometric and Computational Perspectives on Infinite Groups}, G. Baumslag,
D. Epstein, R. Gilman, H. Short, and C. Sims eds., Amer.\ Math.\ Soc., 1995

\bibitem{GH} E. Ghys, P. de la Harpe, eds., Sur les Groupes Hyperboliques 
d'apr\`{e}s Mikhael Gromov, Birkh\"{a}user, 1990. 

\bibitem{Go}
J. Goldstine, Formal languages and their relation to automata: What 
Hopcroft \& Ullman didn't tell us, in {\em Formal Language Theory: 
Perspectives and Open Problems}, R. Book ed., Academic Press, New York 
1980, 109--140.

\bibitem{Gre} R. Gregorac, On generalized free products of finite
extensions
of free groups, J. London Math.\ Soc.\ {\bf 41}, 1966, 662-666.

\bibitem{Gr} M. Gromov, Hyperbolic groups, in{\em Essays in Group Theory},
S.
M. Gersten ed., Springer Verlag 1987, 75--263

\bibitem{HU} J. Hopcroft and J. Ullman, {\em Introduction to Automata
Theory,
Languages, and Computation}, Addison Wesley, 1979.

\bibitem{KPS} A. Karrass, A, Pietrowski, and D. Solitar, An improved 
subgroup theorem for HNN groups with some applicaations, Canadian J.\
Math.\
{\bf 26} 1974, 214-224.

\bibitem{LR} L. P. Lisovik and V. N. Red'ko, Regular events in semigroups, 
Problemy Kibernetiki {\bf 37} 1980, 155-184, 239.

\bibitem{LS} R. C. Lyndon and P. E. Schupp, Combinatorial 
Group Theory, Springer Verlag, 1977.

\bibitem{MS} D. E. Muller and P. E. Schupp, Groups, the theory of ends and
context--free languages, J. Computer and System Sciences {\bf 26} 1983,
295--310.

\bibitem{MS2} D. E. Muller and P. E. Schupp, The theory of ends, pushdown 
autosecond--order logic, Theoretical Computer Sci.\ {\bf 37} 1985, 51-75.

%\bibitem{PS} C. Pittet and L. Saloff--Coste, Amenable groups,
% isoperimetric
% profiles and random walks, Down Under Group Theory, Cossey, Miller,
% Neumann
% and Shapiro, eds., to appear.

\bibitem{SWY}
K. Salomaa, D. Wood, and Sheng Yu, Pumping and pushdown machines, 
Informatique th\'eorique et Applications / Theoretical Informatics and 
Applications, {\bf 28}, 1994, 221--232.

\bibitem{St} J. R. Stallings, Group Theory and Three-Dimensional
Manifolds, Yale Monographs, {\bf 4}, 1971.

\end{thebibliography}

 \end{document}